\documentclass{article}
\usepackage{amssymb}
\usepackage{amsmath}
\usepackage[english]{babel}
\author{E.\,Yu.~Lerner}
\title{Comlexity of prime-dimensional sequences over a finite field}
\date{}
\begin{document}

\maketitle

\begin{abstract}
V.I.~Arnold has recently defined the complexity of a sequence
of~$n$ zeros and ones with the help of the operator of finite
differences. In this paper we describe the results obtained for
almost most complicated sequences of elements of a finite field,
whose dimension~$n$ is a prime number. We prove that with $n\to
\infty$ this property is inherent in almost all sequences, while
the values of multiplicative functions possess this property with
any~$n$ different from the characteristic of the field. We also
describe the prime values of the parameter~$n$ which make the
logarithmic function almost most complicated. All these sequences
reveal a stronger complexity; its algebraic sense is quite clear.
\end{abstract}

\medskip

{\large \bf 1. Main results.}\ In~\cite{ar1}--\cite{ar3},
\cite{ar4} V.I.~Arnold considers the following dynamic system
generated by the finite differentiation operator. Let $x$ be a
closed sequence of $n$ elements from a finite field~$\mathbb{F}_q$
(its $n$th element is followed by the first one). Let $M$ be the
collection of all such sequences ($\# M=q^n$). Let $\Delta: M\to M$ denote the transition from $x$ to the sequence of
differences of neighbor elements of $x$: $x'=\Delta x \Leftrightarrow
x'_i=x_{i+1}-x_i, i=1,\ldots,n$, ($x_{n+1}\equiv x_1$).
The dynamic system $\Delta$ is defined by an oriented graph, whose
vertices are labeled by $x$, $x\in M$. Each vertex $x$ has a
unique outgoing edge (leading to $\Delta x$). The attractors of the
dynamic system $\Delta$ are finite cycles. Each point of the attractor is accessible through a tree of the same
form~(see~\cite{ar1}). V.I.~Arnold has studied the graphs of the dynamic
system $\Delta$ for $q=2$ and $q=3$. Note that the calculations
based on the algorithm described in~\cite{FAN} (see also \cite{Garber})
enable us to obtain these graphs for all values of~$n\le 300$ and $n\le
150$, respectively (see \cite{tabl}, cf. with~\cite{karpen}).

The regular sequences (for example, $(1,\ldots,1)$ or
$(1,0,\ldots,1,0)$) quickly tend to the trivial attractor,
while nonregular ones converge rather slowly to a long cycle. For this reason V.I.~Arnold has defined the complexity of a
sequence of the length~$n$ in terms of the dynamic
system~$\Delta$. The main results of this paper are obtained for prime values
of~$n$. For the mentioned values we prove more accurate statements of
V.I.~Arnold hypotheses on concrete representatives of almost most complicated sequences, as well as the hypotheses on the quota of such
sequences in the entire set~$M$.
In~\cite{ar2} V.I.~Arnold sets up a hypothesis that a random sequence is, as a rule, complicated; its accurate statement and certain related results
are adduced in~\cite{Garber}.

Let us now give several more strict definitions. We understand differential operators as arbitrary linear
operators~$D$ representable in the form $\sum_{i=1}^m d_i \Delta^i$,
$d_i\in \mathbb{F}_q$.
V.I.~Arnold defined the components of the initial vector $x$ with the help of various functions $f$: $x_i=f(i)$. A function $f$ is called {\em almost most complicated}, if $x$ belongs to the attraction
domain of a certain cycle of the maximal period (for a given
map $D$), and the preperiod of the sequence $x,Dx,D^2
x,\ldots$ differs from the maximal one no more than by the unity. With
$D=\Delta$ any almost most complicated function~$f$ is said to be {\em
$\Delta_1$-complicated}.

One can weaken the property of the $\Delta_1$-complexity, neglecting the preperiod of the sequence $x,\Delta x,\Delta^2 x,\ldots$.
A function $f$ such that the corresponding attractor represents a
cycle with the maximal period is called {\em $\Delta_2$-complicated}.
Let $p$ be the {\em characteristic of the field $\mathbb{F}_q$}. According to the results obtained in~\cite{ar1} (see also \cite{Garber}), if~$p$ is an aliquant of $n$, then the notions of $\Delta_1$ and $\Delta_2$-complexity
are equivalent.

V.I.~Arnold has given an example of a complicated function for $q=2$, $n<13$. He
has set up a hypothesis that in the case, when $n+1$ is a certain prime
number~$r$, the following algebraic logarithmic function:
\begin{equation}
\label{lerner1}
f(i)=\left\{
\begin{array}{ll} 0,&\mbox{if $i$ is a quadratic residue modulo $r$,}\\
1,& \mbox{if $i$ is a quadratic nonresidue modulo $r$}
\end{array}
\right.
\end{equation}
is $\Delta_1$-complicated. Unfortunately, certain values of~$n$ make this hypothesis false (see~\cite{Garber}).
The hypothesis on the $\Delta_2$-complexity of this function (which is true for all $n<600$) is still unproved. One can easily verify that with $r$
in the form $4k+3$ this hypothesis is equivalent to the initial one.

A.~Garber has considered the quota of $\Delta_2$-complicated functions among
all possible ones. According to the hypothesis, this quota tends to one as $n\to \infty$. A.~Garber studied
$p$-ary sequences, i.\,e., such that $p=q$. In~\cite{Garber}
he proved that if $n$ has the form $r^k$, where $r$ is a fixed
prime number, then the quota of complicated functions tends to one as
$k\to\infty$.

\medskip

{\large Definition 1.} {\em A function~$f$ is called $D$-complicated,
if it is almost most complicated for \textbf{any}
differential operator.}

\medskip

The algebraic sense of this definition is quite clear (see Lemma~1 and Remark~1). Evidently, this property is more strong than the
$\Delta$-complexity.

\medskip
{\large Theorem~1.} {\em Assume that the dimension $n$ takes on only
prime values. Then the quota of $D$-complicated functions
tends to one as $n\to\infty$.}

\medskip

We prove this theorem constructively, obtaining an explicit formula for the quota of $D$-complicated functions.
\medskip

{\large Theorem~2.}~{\em Let $n=r$ be an odd prime number
different from~$p$. Then the function $f$, which is defined by
formula~(\ref{lerner1}) with $1\le i\le n-1$ and redefined as
$f(n)=0$, is $D$-complicated, if~$p$ is an aliquant of the integer value closest to $n/4$.}

\medskip

In particular, with $q=2$ the function described in Theorem~2 is
$D$-complicated with $n=8k+3$ or $n=8k+5$.

Let $n$ be a prime number. We call a function $f$, mapping $\{1,\ldots,n-1\}$ into
$\mathbb{F}_q$, multiplicative, if $f(ij\ \mbox{mod}\
n)=f(i)f(j)$ for any $i$, $j$ from the definition domain, and $f$ is not
the identical zero. Set~$f(n)=0$. The Legendre symbol~\text{\large
$\left(\frac{i}{n}\right)$} represents an example of such a function.

\medskip
{\large Theorem 3.}\ {\it With\ $n\ne p$ any
multiplicative function is $D$-complicated.}

\medskip

This theorem was first proved in~\cite{FAN}; the proof was elementary but rather intricate. Later, in the journal ``Functional Analysis and Its Applications'' the reviewer of this paper has proposed a very short proof based on more advanced algebraic means, what I sincerely appreciate.
This proof has excited the research described in this paper.

\medskip
{\large \bf 2. An algebra of $D$-complexity.} Let $x$ be a generatrix of a
cyclic group~$C$ of the order~$n$: $x^n=e$, where $e$ is the unit
element of the group. One can identify an arbitrary sequence $f(i)$,
$i=1,\ldots,n$, with an element of the group
algebra\footnote{In \cite{Leng} the corresponding group algebra
is denoted by $\mathbb{F}_q[C]$, here
$\mathbb{F}_q[t]$ stands for a ring of polynomials.}~$A= \mathbb{F}_q C$:
$f=\sum_{i=1}^n f(i) x^i= f(0)+\sum_{i=1}^{n-1} f(i) x^i$,
$f(0)\equiv f(n)$. It is convenient to calculate a product in~$A$,
treating $f$ as a polynomial of the variable~$x$. Then one can calculate the product by the usual multiplication rules for
two polynomials and then reduce the result modulo $x^n-1$. This
simple scheme explicitly defines the isomorphism between $A$ and
$\mathbb{F}_q [t]/(t^n-1)$. The precise representation of~$A$ is also defined by cyclic matrices with $f(0)$ at the main diagonal, with $f(1)$ at the above diagonal, etc.

One can easily see that in terms of the algebra~$A$ a cyclic shift is the multiplication by~$x$; the action of the operator~$\Delta$ is reduced to the multiplication by $x-e$; the action of an arbitrary differential operator is reduced to the multiplication by a fixed element of this algebra which is divisible by $x-e$.

\medskip
{\large Lemma~1.}\ {\it Assume that $n$ is an arbitrary natural
number, $p$ is an aliquant of $n$. A function $f$ is $D$-complicated, if the
corresponding element of the algebra~$A$ is invertible on the subspace
$\sum_{i=1}^n a(i)=0$ of the vector space
$A=(a(1),\ldots,a(n))$.}

\medskip

{\bf Proof of Lemma~1.}\ According to the Chinese remainder theorem,
the algebra of polynomials $\mathbb{F}_q [t]/(t^n-1)$ is representable as the
direct product of algebras $\mathbb{F}_q [t]/(t-1)$ and
$\mathbb{F}_q [t]/\sum_{i=0}^{n-1} t^i$. (Here we use the fact
that $p$ is an aliquant of $n$, consequently, $t=1$ is not a root of the polynomial $\sum_{i=0}^{n-1}t^i$, i.\,e., the latter is not divisible by
$t-1$.) Let $D(t)$ be an
arbitrary polynomial divisible by $(t-1)$; $\tilde f(t)=\sum_{i=0}^{n-1} f(i) t^i$.
In terms of the introduced algebras one can define the complexity of the function~$f$ with the help of the values~$N$ and~$M$, satisfying the equalities
\begin{equation}
\label{system} \left\{\begin{array}{c} D^N(t)\tilde f(t)\!\!\!\mod
\sum_{i=0}^{n-1}t^i=D^M(t) \tilde f(t)\!\!\!\mod \sum_{i=0}^{n-1}t^i,\\
D^N(t) \tilde f(t)\!\!\!\mod (t-1)=D^M(t) \tilde f(t) \!\!\!\mod
(t-1).
\end{array}\right.
\end{equation}
If $F$ is an invertible element of the algebra $\mathbb{F}_q
[t]/\sum_{i=0}^{n-1} t^i$, then the first equality
remains true with the same values of~$N$ and~$M$ for
any function~$f$. This means that both the period and the preperiod of the
sequence $D^N(t)\tilde f(t)\!\!\!\mod
\sum_{i=0}^{n-1}t^i$ are maximum possible. Since the last equality in
system~(\ref{system}) is true with all $N,M\ge 1$, the
function $f$ in this case is $D$-complicated, which was to be proved.

\medskip
{\large Remark~1.}\ {Both the subspace~$S$: $\sum_{i=1}^n a(i)=0$
and the one-dimensional subspace $I$, which is orthogonal to it, are
ideals in~$A$. In more developed algebraic settings, the
invertibility of~$f$ on $S$ means that the projection of~$f$ onto any simple ideal different from $I$ is not zero. In terms of the algebra of
polynomials, this condition is equivalent to the following one: the corresponding polynomial is not divisible by any irreducible polynomial, representing a factor of $\sum_{i=0}^{n-1} t^i$.}
\medskip

{\bf Proof of Theorem~1.} If~$n$ is a prime number, then in accordance with
theorem~2.47 in~\cite{Lidl} the cyclic polynomial
$\sum_{i=0}^{n-1} t^i$ is representable as a product of $(n-1)/d$
different irreducible polynomials of the same degree~$d$, where
$d$ is the order of the number~$q$ in the multiplicative group $F^*_n$.
Consequently, the quota of polynomials mentioned at the end of Remark~1
in the total amount of all polynomials equals $(1-q^{-d})^{(n-1)/d}$. Since, evidently, $q^d\ge n+1$, with $n\to\infty$
we obtain the assertion of the theorem.

\medskip
{\large Remark~2.}\ {Thus, the quota of {\bf cyclic}
matrices which are invertible on the subspace~$S$ tends to the unit,
when $n$ tends to infinity (taking on only prime values). It is interesting
that the quota of {\bf various} matrices
invertible on this subspace tends to a certain value from $(0, 1)$ (see, for example,
\cite{OBLUP}).}

\medskip
{\bf Proof of Theorem~2.}\ Let $\zeta$ be a primitive $n$th
root of the unit in the corresponding algebraic
extension of the field~$\mathbb{F}_q$ (see \cite[Chap.~8, \S~3]{Leng}).
One can easily verify~(see~\cite{Sachkov}) that eigenvalues of a
cyclic matrix with the first row $(f(0),\ldots,f(n-1))$, are
$$
\lambda_m=\sum_{i=1}^n f(i)\zeta^{im},\ m=1,\ldots,n,
$$
and the corresponding eigenvectors take the form
$$
(1,\zeta^m,\ldots,\zeta^{m(n-1)}),\quad m=1,\ldots,n.
$$
The first~$n-1$ eigenvectors generate the space~$S$. In accordance with
Lemma~1 the function~$f$ is $D$-complicated, if the product of the
corresponding eigenvalues differs from zero.

Let $g_m=\sum_{j=1}^n \text{\large $\left(\frac{j}{n}\right)$}
 \zeta^{jm}$ be the Gauss sums. Due to the multiplicative property of the Legendre symbol we have $g_m=\text{\large $\left(\frac{m}{n}\right)$}
 g_1$. It is well known that (see \cite[Chap.~8, \S~3]{Leng})
 $g_1^2=\text{\large $\left(\frac{-1}{n}\right) n$}$, and, evidently,
$\sum_{i=1}^{n-1}\zeta^{im}=-1$ with $m=1,\ldots,n-1$. Using these
equalities, we obtain that the function~$f$ mentioned in the assumption of the theorem satisfies the relation
$$
\prod_{m=1}^{n-1} \lambda_m= \left\{
\begin{array}{ll} k^{(n-1)/2},&\text{if $n=4k+1$,}\\
(k+1)^{(n-1)/2},& \text{if $n=4k+3$.}
\end{array}
\right.
$$
The theorem is proved.

\medskip
{\bf Proof of Theorem~3.}\ Without loss of generality, we assume
that the function $f$ differs from the $\delta$-function
$(1,0,\ldots,0)$, and, consequently, $n>2$. An arbitrary
nontrivial automorphism of a cyclic group defined by the formula
$x\to x^k$, $k\!\!\!\mod n\ne 0,1$, transitively represents all
elements of this group different from $e$, as well as all simple
ideals of the algebra~$A$ which differ from~$I$. This transform maps the
element~$f=\sum_{i=1}^{n-1} f(i) x^i$ which corresponds to the
multiplicative function into $\sum_{i=1}^{n-1} f(ik) x^i=f(k) f$.
Therefore, if $f(k)\ne 0$ for certain~$k>1$, then either the
projections of~$f$ onto all simple ideals different from~$I$ are not zeros, or the projection of $f$ on $S$ equals zero. The latter means that
$f=c \sum_{i=0}^{n-1} x^i$, what contradicts the conditions~$f(0)=0$,
$f\not\equiv 0$. The theorem is proved.

\end{document}